\magnification=1200
\overfullrule=0pt
\centerline {\bf Revisiting a theorem on multifunctions of one real variable}\par
\bigskip
\bigskip\centerline {BIAGIO RICCERI}\par
\bigskip
\bigskip
\centerline {\it Dedicated to Professor Simeon Reich, with esteem, on his 65th birthday}\par
\bigskip
\bigskip
{\bf Abstract:} In this paper, 
 we intend to revisit Theorem 2 of [3]  formulating it in a way that, weakening the hypotheses 
and, at the same time,  highlighting the richer
conclusion allowed by the proof, it can potentially be applicable to a broader range of different situations.
Samples of such applications are also given.\par
 \par
\bigskip
\bigskip
\bigskip
\bigskip
Some years ago, we established a certain theorem ([3], Theorem 2) on a class of multifunctions depending on a real
variable whose formulation was heavily conditioned by the application of it to minimax theory which just  was the core of
[3].\par
\smallskip
In the present paper, we intend to revisit that result formulating it in a way that, weakening the hypotheses 
and, at the same time,  highlighting the
 richer
conclusion allowed by the proof, it can potentially be applicable to a broader range of different situations.\par
\smallskip
So, after establishing the main result (Theorem 1), we give a sample of application of it (Theorem 3) that
cannot be deduced by Theorem 2 of [3]. In turn, we highlight a series of consequences of Theorem 3 essentially
dealing with the existence of some kind of "singular" points for functions of the type $f+\lambda g$, with $\lambda\in {\bf R}$.\par
\smallskip
In the sequel, the term "interval" means a non-empty connected subset of ${\bf R}$ with more than one point.\par
\smallskip
  For 
a multifunction $F:I\to 2^X$, as usual, for $A\subseteq I$ and $B\subseteq X$, we set
$$F(A)=\cup_{x\in A}F(x)$$
and
$$F^{-}(B)=\{\lambda\in I : F(\lambda)\cap B\neq\emptyset\}\ .$$
When $I$ is an interval,  $F$ 
  is said to be non-decreasing (resp.
non-increasing) with respect to the inclusion if $F(\lambda)\subseteq F(\mu)$ 
(resp. $F(\mu)\subseteq F(\lambda))$ for all $\lambda, \mu\in I$, with $\lambda<\mu$.\par
\smallskip
We start by proving the following\par
\medskip
PROPOSITION 1. - {\it Let $X, Y$ be two non-empty sets,  $D\subseteq Y$, 
$F:X\to 2^Y$ a
multifunction such that $F(x)\cap D\neq\emptyset$ for all $x\in X$. Assume also that
there exist $y_0\in Y$ and a topology on $X\setminus F^-(y_0)$
 such that $X\setminus F^{-}(y_0)$ is sequentially
compact (resp. compact) and $X\setminus F^{-}(\{y,y_0\})$ is sequentially
closed (resp. closed) in $X\setminus F^{-}(y_0)$ for all $y\in D$.\par
Then,  for every
 non-decreasing sequence $\{Y_n\}$ of subsets of $Y$,
with $\cup_{n\in {\bf N}}Y_n=Y$, 
there exists $\tilde n\in {\bf N}$ such that
$F(x)\cap Y_{\tilde n}\neq\emptyset$ for all $x\in X$.}\par
\smallskip
PROOF. Let $\{Y_n\}$ be a
 non-decreasing sequence of subsets of $Y$,
with $\cup_{n\in {\bf N}}Y_n=Y$. 
 Fix $\nu\in {\bf N}$
so that $y_0\in Y_{\nu}$.
Arguing by contradiction, assume
that, for each $n\in {\bf N}$, there exists $x_n\in X$ such that
$$F(x_n)\cap Y_n=\emptyset\ .\eqno{(1)}$$
First, consider the ''sequentially compact, sequentially closed" case.
 Hence, for each $n\geq\nu$, one
has $y_0\not\in F(x_n)$, that is $x_n\in X\setminus F^{-}(y_0)$.
So, there exists a subsequence $\{x_{n_k}\}$ converging to 
a point $x^*\in X\setminus F^{-}(y_0)$. Now, fix $y^*\in F(x^*)\cap D$
and $h\geq\nu$ such that $y^*\in Y_h$. By assumption, $F^{-}(y^*)\cap
(X\setminus F^{-}(y_0))$ is sequentially open in $X\setminus F^{-}(y_0)$, and
 hence $x_{n_k}\in F^{-}(y^*)$ for all $k$ large enough. Then, if we choose
$k$ so that $n_k\geq h$, we have $y^*\in F(x_{n_k}\cap Y_{n_k})$, against
$(1)$.  Now, consider the "compact, closed" case. Let $A\subseteq D$ be a finite set.
Fix $p\geq \nu$ so that $A\subseteq Y_p$. Hence, in view of $(1)$, we have
$$X\setminus  F^{-}(A\cup \{y_0\})\neq 0 \ .$$
In other words, the family $\{(X\setminus F^{-}(y))\cap (X\setminus F^{-}(y_0))\}_{y\in D}$
has the finite intersection property. But then, since each member of this
family is closed in $X\setminus F^-(y_0)$ which is compact, we have
$$X\setminus F^{-}(D\cup \{y_0\})\neq 0 \ .$$
This is against the assumption that $F^-(D)=X$, and the proof is complete.
\hfill $\bigtriangleup$\par
\medskip
Our main result is as follows.\par
\medskip
THEOREM 1. - {\it Let $X$ be a non-empty set, $I\subseteq {\bf R}$ an interval and
$F:I\to 2^X$ a multifunction satisfying the following conditions:\par
\noindent
$(i)$\hskip 5pt there exist $\lambda_0\in I$, with $F(\lambda_0)\neq\emptyset$, and
a topology on $F(\lambda_0)$ such that 
$F(\lambda_0)$ is sequentially compact (resp. compact);\par
\noindent
$(ii)$\hskip 5pt the set 
$$\{\lambda\in I : 
F(\lambda)\cap F(\lambda_0)\hskip 3pt is\hskip 3pt
 sequentially\hskip 3pt closed\hskip 3pt (resp.\hskip 3pt closed)\hskip 3pt in\hskip 3pt F(\lambda_0)\}$$
 is dense in $I$\ ;\par
\noindent
$(iii)$\hskip 5pt for each $x\in X$, the set $I\setminus F^{-}(x)$ is an
 interval open in $I$\ .\par
Under such hypotheses, there exists a compact interval
$[a^*,b^*]\subseteq I$  such that either
$( F(a^*)\cap F(\lambda_0))\setminus F(]a^*,b^*[)\neq \emptyset$ and $F_{|]a^*,b^*[}$
is non-decreasing with respect to the inclusion, or 
$( F(b^*)\cap F(\lambda_0))\setminus F(]a^*,b^*[)\neq \emptyset$ and $F_{|]a^*,b^*[}$  
is non-increasing with respect to the inclusion. 
In particular, the first (resp. second)
occurrence is true when $\lambda_0=\inf I$ (resp. $\lambda_0=\sup I$)\ .}\par
\smallskip
 PROOF. For each $x\in X$, put
$$\Phi(x)=I\setminus F^{-}(x)\ .$$
Clearly
$$\Phi^{-}(\lambda)=X\setminus F(\lambda)$$
for all $\lambda\in I$.  In view of Proposition 1,
there exists a compact interval $[a,b]\subseteq I$, with $\lambda_0\in
[a,b]$, such that
$$\Phi(x)\cap [a,b]\neq \emptyset$$
for all $x\in X$. Therefore, each set $\Phi(x)\cap [a,b]$ is an
interval open in $[a,b]$. For each $x\in X$, put
$$\alpha(x)=\inf(\Phi(x)\cap [a,b])$$
and
$$\beta(x)=\sup(\Phi(x)\cap [a,b])\ .$$
Clearly, for each $x_0\in F(\lambda_0)$ and each $r\in ]\alpha(x_0),\beta(x_0)[\cap D$
one has
$$x_0\in \Phi^{-}(r)$$
and
$$\alpha(x)<r<\beta(x)$$
for all $x\in \Phi^{-}(r)$. Since, by assumption,  $\Phi^{-}(r)\cap F(\lambda_0)$ is sequentially open  (resp. open)
in $F(\lambda_0)$
and
$D$ is dense in $I$, we then infer that $\alpha_{|F(\lambda_0)}$ is sequentially upper semicontinuous (resp.
upper semicontinuous)
at $x_0$, while $\beta_{|F(\lambda_0)}$ is sequentially lower semicontinuous (resp. lower semicontinuous) at $x_0$. Now,
suppose that $\lambda_0\in ]a,b[$. Observe that
$$F(\lambda_0)=\alpha^{-1}([\lambda_0,+\infty[)\cup \beta^{-1}(]-\infty,\lambda_0]) \ .\eqno{(2)}$$
Since $F(\lambda_0)\neq\emptyset$, we have either
 $\alpha^{-1}([\lambda_0,+\infty[)\neq\emptyset$ or 
$\beta^{-1}(]-\infty,\lambda_0])\neq\emptyset$.  First,
assume that   $\alpha^{-1}([\lambda_0,+\infty[)\neq\emptyset$. Then, since
$F(\lambda_0)$ is sequentially compact (resp. compact) and $\alpha_{|F(\lambda_0)}$ is sequentially upper
semicontinuous (resp. upper semicontinuous), in view of $(2)$, there is $x^*\in F(\lambda_0)$ such that $\alpha(x^*)=\sup_X\alpha$.
Since $\alpha(x^*)\geq \lambda_0$, we have $\alpha(x^*)\in ]a,b[$. This implies,
in particular, that $\alpha(x^*)$ does not belong to $\Phi(x^*)\cap [a,b]$, since this set
is open in $[a,b]$. As a consequence, we have $x^*\in F(\alpha(x^*))$.  Now, fix
$\lambda, \mu\in ]\alpha(x^*),\beta(x^*)[$, with $\lambda<\mu$.  Clearly,  $\mu\not\in F^{-}(x^*)$ and
hence $x^*\not\in F(\mu)$. Next, for each $x\in \Phi^{-}(\mu)$,
we have
$$\alpha(x)\leq \alpha(x^*)<\lambda<\mu\leq\beta(x)\ .$$
Hence, $\lambda\not\in F^{-}(x)$ that is $x\in \Phi^{-}(\lambda)$. Therefore,
we have
$$x^*\in F(\alpha(x^*))\setminus F(]\alpha(x^*),\beta(x^*)[)$$
as well as
$$\Phi^{-}(\mu)\subseteq \Phi^{-}(\lambda)$$
that is
$$F(\lambda)\subseteq F(\mu)\ .$$
So, in the current case, the conclusion is satisfied taking $a^*=\alpha(x^*)$ and $b^*=\beta(x^*)$.
Now, assume that $\beta^{-1}(]-\infty,\lambda_0])$ is non-empty. This time, due to
the sequential lower semicontinuity (resp. lower semicontinuity) of $\beta_{|F(\lambda_0)}$, there exists $\hat x\in X$ such
that $\beta(\hat x)=\inf_X\beta$. As before, one realizes that $\hat x\in F(\beta(\hat x))$.
Fix $\lambda, \mu\in ]\alpha(\hat x),\beta(\hat x)[$ with $\lambda<\mu$. Clearly,
$\hat x\not\in F(\lambda)$. For each $x\in \Phi^{-}(\lambda)$, we have
$$\alpha(x)\leq \lambda<\mu<\beta(\hat x)\leq \beta(x)$$
and so $x\in \Phi^{-}(\mu)$. Therefore, we have
$$\hat x\in F(\alpha(\hat x))\setminus F(]\alpha(\hat x),\beta(\hat x)[)$$
as well as
$$\Phi^{-}(\lambda)\subseteq \Phi^{-}(\mu)$$
that is
$$F(\mu)\subseteq F(\lambda)\ .$$
So, in this case, the conclusion is satisfied taking $a^*=\alpha(\hat x)$ and $b^*=\beta(\hat x)$.
Now, assume  $\lambda_0=a$ (in particular, this occurs when
$\lambda_0=\inf I$).  
If $\sup_X\alpha>a$, then since
$$\alpha^{-1}([a,+\infty[)\subseteq F(a)\ ,$$
and $F(a)$ is sequentially compact (resp. compact), $\alpha$ attains its supremum (larger than $a$), and so
we are exactly in the  first sub-case considered when $\lambda_0\in ]a,b[$. If $\sup_X\alpha=a$,  
we still reach
the conclusion, as in the first sub-case considered when $\lambda_0\in ]a,b[$, taking  $a^*=a$ and
$b^*=\beta(x^*)$, where $x^*$ is any point in $F(a)$.
 In doing so, notice simply that $x^*\in F(\alpha(x^*))$ as $\alpha(x^*)=a$.
Finally, let $\lambda_0=b$ ( in particular, this occurs when $\sup I=b$). 
If $\inf_X\beta<b$, then since
$$\beta^{-1}(]-\infty,b])\subseteq F(b)\ ,$$
and $F(b)$ is sequentially compact (resp. compact), $\beta$ attains its supremum (smaller than $b$), and so
we are exactly in the second sub-case considered when $\lambda_0\in ]a,b[$. 
If $\inf_X\beta=b$, we still reach
the conclusion, as in the second sub-case considered when $\lambda_0\in ]a,b[$, taking  $a^*=\alpha(\hat x)$ and
$b^*=b$, where $\hat x$ is any point in $F(b)$.  The proof is complete.\hfill $\bigtriangleup$\par
\medskip
We now give a purely set-theoretical reformulation of Theorem 1 in the "compact, closed" case. We first need
the following definition.\par
\medskip
DEFINITION 1. - Let $Y$ be a non-empty set and ${\cal F}$ a family of
subsets of $Y$. We say that ${\cal F}$ has the compactness-like property if
every subfamily of ${\cal F}$ satisfying the finite intersection property has
a non-empty intersection. \par
\medskip
We have the following characterization which is due to C. Costantini ([1]):\par
\medskip
PROPOSITION 2. - {\it Let $Y$ be a non-empty set, let
${\cal F}$ be a family of subsets of $Y$ and let $\tau$ be the topology
on $Y$ generated by the family $\{Y\setminus C\}_{C\in {\cal F}}$.\par
Then, the following assertions are equivalent:\par
\noindent
$(a)$\hskip 5pt Each member of ${\cal F}$ is $\tau$-compact.\par
\noindent
$(b)$\hskip 5pt The family ${\cal F}$ has the compactness-like property.\par
\noindent
$(c)$\hskip 5pt The space $Y$ is $\tau$-compact.}\par
\medskip
Here is the reformulation of Theorem 1.\par
\medskip
THEOREM 2. - {\it Let $X$ be a non-empty set, $I\subseteq {\bf R}$ a
 interval and $F:I\to 2^X$ a multifunction such that,
for each $x\in X$, the set $X\setminus F^-(x)$ is an interval
open in $I$. Moreover, assume that, for some $\lambda_0\in I$, with
$F(\lambda_0)\neq\emptyset$, and some
set $D\subseteq I$ dense in $I$, the family $\{F(\lambda)\cap F(\lambda_0)\}_{\lambda\in D}$
has the compactness-like property.\par
Then, the same conclusion as that of Theorem 1 holds.}\par
\smallskip
PROOF. In view of Proposition 2, if we consider the topology on $F(\lambda_0)$
generated by the family $\{F(\lambda_0)\setminus F(\lambda)\}_{y\in D}$, all the
assumptions of Theorem 1 (for the "compact, closed" case) are satisfied, and
the conclusion follows.\hfill $\bigtriangleup$.\par
\medskip
Now, we are going to present an application of Theorem 1. \par
\smallskip
In the sequel, $X$ is a non-empty set, $Y$ is a real Hausdorff locally convex
topological vector space, $C$ is a closed subset of $Y$ such that $Y\setminus C$ is convex,
$I\subseteq {\bf R}$ is an interval containing $0$ and $f, g$ are two functions from $X$ into
$Y$. The symbol $\partial$ stands for boundary.\par
\smallskip
Here is the above mentioned application.
\medskip
THEOREM 3. - {\it 
Assume that the following conditions are satisfied: \par
\noindent
$(a_1)$\hskip 5pt the set $f^{-1}(C)$
 is non-empty and the set $\{(f(x),g(x)) : x\in f^{-1}(C)\}$ is compact in $Y\times Y$\ ;\par
\noindent
$(a_2)$\hskip 5pt for each $x\in X$,  
there exists $\lambda\in I$ such that
$$f(x)+\lambda g(x)\in Y\setminus C\ .$$
Then, there exist a  compact interval $[a^*,b^*]\subseteq I$ and
a point $x^*\in f^{-1}(C)$ satisfying
$$f(x^*)+\lambda g(x^*)\in Y\setminus C$$
for all $\lambda\in ]a^*,b^*[$, such that, if we put
$$V=\bigcup_{\lambda\in ]a^*,b^*[}\{x\in X : f(x)+\lambda g(x)\in Y\setminus C\}\ ,$$ 
 at least one of the following assertions holds:\par
\noindent
$(p_1)$\hskip 5pt $f(x^*)+a^*g(x^*)\in \partial C$ 
and 
$$(f+a^*g)(V)\cap C\subseteq\partial (f+a^*g)(V)\cap \partial C\ ;$$
$(p_2)$\hskip 5pt $f(x^*)+b^*g(x^*)\in \partial C$ 
and 
$$(f+b^*g)(V)\cap C\subseteq\partial (f+b^*g)(V)\cap \partial C\ .$$
In particular, $(p_1)$ (resp. $(p_2)$)
holds when $0=\inf I$ (resp. $0=\sup I$)\ .}\par
\smallskip
PROOF. Consider the multifunction
 $F:I\to 2^X$  defined by
$$F(\lambda)=\{x\in X : f(x)+\lambda g(x)\in C\}$$
for all $\lambda\in I$. Observe that, taking $\lambda_0=0$, $F$
satisfies the assumptions of Theorem 1. Indeed,
if we consider on $f^{-1}(C)$ the weakest topology for which both $f$ and $g$
are continuous in $f^{-1}(C)$, then, in view of $(a_1)$, $f^{-1}(C)$ turns out to
be compact in that topology. So, $(i)$ is satisfied. Since $Y$ carries a vector topology,
for each $\lambda\in {\bf R}$, the function $f+\lambda g$ is continuous in $f^{-1}(C)$,
and so also $(ii)$ is satisfied since $C$ is closed. Finally, 
for each $x\in X$, we have
$$I\setminus F^-(x)=\{\lambda\in I : f(x)+\lambda g(x)\in Y\setminus C\}$$
which is an interval open $I$, in view of $(a_2)$ and of the fact that
$Y\setminus C$ is open and convex. 
Therefore, Theorem 1 ensures the existence of a 
compact interval $[a^*,b^*]\subseteq I$ such that either 
$( F(a^*)\cap F(0))\setminus F(]a^*,b^*[)\neq \emptyset$ and $F_{|]a^*,b^*[}$
is non-decreasing with respect to the inclusion, or 
$( F(b^*)\cap F(0))\setminus F(]a^*,b^*[)\neq \emptyset$ and $F_{|]a^*,b^*[}$  
is non-increasing with respect to the inclusion. 
Assume, for instance, that
$( F(a^*)\cap F(0))\setminus F(]a^*,b^*[)\neq \emptyset$ and $F_{|]a^*,b^*[}$
is non-decreasing with respect to the inclusion. Pick $x^*\in ( F(a^*)\cap F(0))\setminus F(]a^*,b^*[)$.
So, $x^*\in f^{-1}(C)$, $f(x^*)+a^*g(x^*)\in C$ and $f(x^*)+\lambda g(x^*)\in Y\setminus C$ for all $\lambda\in ]a^*,b^*[$.
Now, let us show that $(f+a^*g)(V)\cap C\subseteq \partial C$. So, let
$x\in V$ be such that $f(x)+a^*g(x)\in C$. Fix $\lambda\in ]a^*,b^*[$ such that
$f(x)+\lambda g(x)\in Y\setminus C$.   Arguing by contradiction,
suppose that $f(x)+a^*g(x)\in$ int$(C)$. Then, we could find $\delta\in
]a^*,\lambda[$ so that $f(x)+\delta g(x)\in \hbox {\rm int}(C)$. But this contradicts the fact
that $f(x)+\delta g(x)\subseteq Y\setminus C$ as $F(\delta)\subseteq F(\lambda)$. 
Now, let us show that $(f+a^*g)(V)\cap C\subseteq\partial (f+a^*g)(V)$. So, let $z\in (f+a^*g)(V)\cap C$.
Arguing by contradiction again, assume that $z\in$ int$((f+a^*g)(V))$. 
Since $Y\setminus C$ is open and convex and $z\in \partial (Y\setminus C)$, there exists
 $\varphi\in Y^*\setminus \{0\}$ such that $\varphi(z)\leq
\varphi(u)$ for all $u\in Y\setminus C$. Hence, the set $\varphi^{-1}(]-\infty, \varphi(z)[)$
is an open set contained in $C$ which meets int$((f+a^*g)(V))$ since $\varphi$ has no
local minima being linear, and this is impossible for what seen above. The proof is complete.
\hfill $\bigtriangleup$\par
\medskip
REMARK 1. -  Let us recall that a function $h:X\to Y$ between topological spaces is said to be
open at $x_0\in X$ if there exists a fundamental system ${\cal V}$ of neighbourhoods of $x_0$
such that, for each $V\in {\cal V}$, the set $h(V)$ is a neighbourhood of $h(x_0)$. Now, in connection
with Theorem 3, if $\tau$ is any topology on $X$ containing the family
$$\{\{x\in X : f(x)+\lambda g(x)\in Y\setminus C\}\}_{\lambda\in I}\ ,$$
it follows that at least one of the functions $f+a^*g$, $f+b^*g$ is not $\tau$-open at the point
$x^*$.\par
\medskip
Among the corollaries of Theorem 3, it is worth noticing the following\par
\medskip
THEOREM 4. - {\it 
 Let $\varphi\in Y^*\setminus
\{0\}$ and $r\in {\bf R}$ be such that the set
$$K:=\{x\in X : \varphi(f(x))\leq r\}$$
is non-empty. Assume also that the set $\{(f(x),g(x)) : x\in K\}$ is compact
in $Y\times Y$ and that $g(K)\cap
\varphi^{-1}(0)=\emptyset$.\par
Then, there exist a  compact interval $[a^*,b^*]\subseteq {\bf R}$ and
a point $x^*\in K$ satisfying
$$\varphi(f(x^*)+\lambda g(x^*))>r$$
for all $\lambda\in ]a^*,b^*[$, such that, if we put
$$V=\bigcup_{\lambda\in ]a^*,b^*[}\{x\in X : \varphi(f(x)+\lambda g(x))>r\}\ ,$$ 
 at least one of the following assertions holds:\par
\noindent
$(q_1)$\hskip 5pt $\varphi(f(x^*)+a^*g(x^*))=r$ 
and 
$$(f+a^*g)(V)\cap \varphi^{-1}(]-\infty,r])\subseteq
\partial (f+a^*g)(V)\cap \varphi^{-1}(r)\ ;$$
$(q_2)$\hskip 5pt $\varphi(f(x^*)+b^*g(x^*))=r$ 
and 
$$(f+b^*g)(V)\cap \varphi^{-1}(]-\infty,r])\subseteq
\partial (f+b^*g)(V)\cap \varphi^{-1}(r)\ .$$}
\smallskip
PROOF. It is enough to apply Theorem 3 taking $I={\bf R}$ and
 $C=\varphi^{-1}(]-\infty,r])$, so
that $\partial C=\varphi^{-1}(r)$.\par
\medskip
If we apply Theorem 3 jointly with [2], we obtain:\par
\medskip
THEOREM 5. - {\it Let $Y$ be a finite-dimensional Banach space and
let $\psi:Y\to Y$ be a continuous
function such that $\psi^{-1}(C)$ is non-empty and compact. Assume also that,
for each $x\in X$, there exists $\lambda\in {\bf R}$ such that
$$\psi(x)+\lambda x\in X\setminus C\ .$$
Then, there exist $x^*\in \psi^{-1}(C)$ and $\mu^*\in {\bf R}$ such that
$$\psi(x^*)+(1+\mu^*)x^*\in\partial C$$
and
$$\sup_{\|x-x^*\|\leq r}\|\psi(x)-\psi(x^*)+\mu^*(x-x^*)\|\geq r$$
for each $r>0$ small enough.}\par
\smallskip
PROOF. Apply Theorem 3 taking $X=Y$, $I={\bf R}$, $f=\psi$ and $g=$id.  Then, there exist $x^*\in \psi^{-1}(C)$,
$\lambda^*\in {\bf R}$ and a neighbourhood $V$ of $x^*$  such that 
$$\psi(x^*)+\lambda^*x^*\in \partial(\psi+\lambda^*\hbox {\rm id})(V)\cap \partial C\ .$$
 Now, from the proof of Theorem 1
of [2], we know that if $r>0$ is such that
$$\sup_{\|x-x^*\|\leq r}\|\psi(x)-\psi(x^*)+(\lambda^*-1)(x-x^*)\|<r\ ,$$
then, for some $r_0>0$ and for every $y\in X$ satisfying $\|\psi(x^*)+\lambda^* x^*-y\|<r_0$, there
exists $x\in X$, with $\|x-x^*\|\leq r$, such that 
$$\psi(x)+\lambda^* x=y\ .$$
As a consequence, for every $r>0$ for which the closed ball centered at $x^*$, of radius $r$, is contained
in $V$, we have
$$\sup_{\|x-x^*\|\leq r}\|\psi(x)-\psi(x^*)+(\lambda^*-1)(x-x^*)\|\geq r\ .$$
Now, the conclusion follows taking $\mu^*=\lambda^*-1$.\hfill $\bigtriangleup$
\medskip
The last consequence of Theorem 3 that we point out is as follows:\par
\medskip
THEOREM 6. - {\it Let $X$ be an open set in a real Banach space, let $Y$ be a Banach space
 and let $f, g$ be continuously Fr\'echet differentiable. Assume that $(a_1)$, $(a_2)$ hold.\par
Then, there exist $x^*\in X$ and $\lambda^*\in I$ such that
$$f(x^*)+\lambda^*g(x^*)\in \partial C$$
and the continuous linear operator $f'(x^*)+\lambda^*g'(x^*)$ is not invertible.}\par
\smallskip
PROOF. By Theorem 3, there exists $x^*\in X$ and $\lambda^*\in I$ such that
$$f(x^*)+\lambda^*g(x^*)\in \partial C$$
and $f+\lambda^*g$ is not open at $x^*$ (Remark 1). This just implies that
$f'(x^*)+\lambda^*g'(x^*)$ is not invertible, since, otherwise, $f+\lambda^*g$
would be a local homeomorphism at $x^*$ by the inverse function theorem.\hfill
$\bigtriangleup$
\bigskip
\bigskip
\bigskip
\bigskip
\centerline {\bf References}\par
\bigskip
\bigskip
\noindent
[1]\hskip 5pt C. COSTANTINI, Personal communication.\par
\smallskip
\noindent
[2]\hskip 5pt M. REICHBACH, {\it
 Fixed points and openness},  
Proc. Amer. Math. Soc., {\bf 12} (1961), 734-736.\par
\smallskip
\noindent
[3]\hskip 5pt B. RICCERI, {\it A further improvement of a minimax theorem of
Borenshtein and Shul'man}, J. Nonlinear Convex Anal., {\bf 2} (2001),
279-283.\par
\bigskip
\bigskip
\bigskip
\bigskip
Department of Mathematics\par
University of Catania\par
Viale A. Doria 6\par
95125 Catania\par
Italy\par
{\it e-mail address}: ricceri@dmi.unict.it

\bye